
\catcode`\@=11
\def\undefine#1{\let#1\undefined}
\def\newsymbol#1#2#3#4#5{\let\next@\relax
 \ifnum#2=\@ne\let\next@\msafam@\else
 \ifnum#2=\tw@\let\next@\msbfam@\fi\fi
 \mathchardef#1="#3\next@#4#5}
\def\mathhexbox@#1#2#3{\relax
 \ifmmode\mathpalette{}{\m@th\mathchar"#1#2#3}%
 \else\leavevmode\hbox{$\m@th\mathchar"#1#2#3$}\fi}
\def\hexnumber@#1{\ifcase#1 0\or 1\or 2\or 3\or 4\or 5\or 6\or 7\or 8\or
 9\or A\or B\or C\or D\or E\or F\fi}

\newdimen\ex@
\ex@.2326ex
\def\varinjlim{\mathop{\vtop{\ialign{##\crcr
 \hfil\rm lim\hfil\crcr\noalign{\nointerlineskip}\rightarrowfill\crcr
 \noalign{\nointerlineskip\kern-\ex@}\crcr}}}}
\def\varprojlim{\mathop{\vtop{\ialign{##\crcr
 \hfil\rm lim\hfil\crcr\noalign{\nointerlineskip}\leftarrowfill\crcr
 \noalign{\nointerlineskip\kern-\ex@}\crcr}}}}
\def\varliminf{\mathop{\underline{\vrule height\z@ depth.2exwidth\z@
 \hbox{\rm lim}}}}

\font\tenmsa=msam10 \font\sevenmsa=msam7 \font\fivemsa=msam5
\newfam\msafam
\textfont\msafam=\tenmsa \scriptfont\msafam=\sevenmsa
\scriptscriptfont\msafam=\fivemsa
\edef\msafam@{\hexnumber@\msafam} \mathchardef\dabar@"0\msafam@39
\def\dashrightarrow{\mathrel{\dabar@\dabar@\mathchar"0\msafam@4B}}
\def\dashleftarrow{\mathrel{\mathchar"0\msafam@4C\dabar@\dabar@}}

\font\tenmsb=msbm10 \font\sevenmsb=msbm7 \font\fivemsb=msbm5
\newfam\msbfam
\textfont\msbfam=\tenmsb \scriptfont\msbfam=\sevenmsb
\scriptscriptfont\msbfam=\fivemsb
\edef\msbfam@{\hexnumber@\msbfam}
\def\Bbb#1{{\fam\msbfam\relax#1}}
\def\widehat#1{\setbox\z@\hbox{$\m@th#1$}%
 \ifdim\wd\z@>\tw@ em\mathaccent"0\msbfam@5B{#1}%
 \else\mathaccent"0362{#1}\fi}
\font\teneufm=eufm10 \font\seveneufm=eufm7 \font\fiveeufm=eufm5
\newfam\eufmfam
\textfont\eufmfam=\teneufm \scriptfont\eufmfam=\seveneufm
\scriptscriptfont\eufmfam=\fiveeufm

\newsymbol\boxtimes 1202
\newsymbol\rtimes 226F

\catcode`\@=12

\magnification=\magstep1 \font\title = cmr10 scaled \magstep2
 
\font\author = cmcsc10 \font\addr = cmti8 \font\byabs = cmr8
  \font\email = cmtt8

\parindent=1em
\baselineskip 15pt
\vsize=18.5 cm

\newcount\refcount
\newcount\seccount
\newcount\sscount
\newcount\eqcount
\newcount\boxcount
\newcount\testcount
\newcount\bibcount
\boxcount = 128 \seccount = -1

\def\sec#1{\advance\seccount by 1\bigskip\goodbreak\noindent
    {\bf\number\seccount.\ #1}\medskip \sscount = 0\eqcount = 0}
\def\proc#1#2{\advance\sscount by 1\eqcount = 0
    \medskip\goodbreak\noindent{\author #1}
    {\tenrm{\number\sscount}}:\ \ {\it #2}}
\def\nproc#1#2#3{\advance\sscount by 1\eqcount = 0\global
    \edef#1{#2\ \number\sscount}
    \medskip\goodbreak\noindent{\author #2}
    {\tenrm{\number\sscount}}:\ \ {\it #3}}
\def\lemma#1{\advance\sscount by 1\eqcount = 0
    \medskip\goodbreak\noindent{\author Lemma}
    {\tenrm{\number\sscount}}:\ \ {\it #1}}
\def\nlemma#1#2{\advance\sscount by 1\eqcount = 0\global
    \edef#1{Lemma\ \number\sscount}
    \medskip\goodbreak\noindent{\author Lemma}
    {\tenrm{\number\sscount}}:\ \ {\it #2}}
\def\proof{\medskip\noindent{\it Proof:\ \ }}

\def\eql#1{\global\advance\eqcount by 1\global
    \edef#1{(\number\sscount.\number\eqcount)}\leqno{#1}}
\def\ref#1#2{\advance\refcount by 1\global
    \edef#1{[\number\refcount]}\setbox\boxcount=
    \vbox{\item{[\number\refcount]}#2}\advance\boxcount by 1}
\def\biblio{{\frenchspacing
    \bigskip\goodbreak\centerline{\it References}\medskip
    \bibcount = 128\loop\ifnum\testcount < \refcount
    \goodbreak\advance\testcount by 1\box\bibcount
    \advance\bibcount by 1\vskip 4pt\repeat\medskip}}
\def\tightmatrix#1{\null\,\vcenter{\normalbaselines\mathsurround=0pt
    \ialign{\hfil$##$\hfil&&\ \hfil$##$\hfil\crcr
    \mathstrut\crcr\noalign{\kern-\baselineskip}
    #1\crcr\mathstrut\crcr\noalign{\kern-\baselineskip}}}\,}

\def\emph{\bf}
\def\colon{{:}\;}
\def\|{|\;}


\def\C{{\Bbb C}}
\def\Z{{\Bbb Z}}

\def\F{{\Bbb F}}

\def\N{{\Bbb N}}

\def\G{{\Bbb G}}

\def\Aut{{\rm Aut}}
\def\Out{{\rm Out}}

\def\End{{\rm End}}
\def\Ind{{\rm Ind}}
\def\Res{{\rm Res}}

\def\PGL{{\rm PGL}}
\def\GL{{\rm GL}}

\def\SU{{\rm SU}}
\def\SO{{\rm SO}}
\def\SL{{\rm SL}}

\def\PSL{{\rm PSL}}
\def\PSU{{\rm PSU}}

\def\qed{\hfill\hbox{$\sqcup$\llap{$\sqcap$}}\medskip}

\def\m#1#2#3#4{\left(\matrix{#1&#2\cr #3&#4\cr}\right)}
\def\F{{\bf F}}
\def\Z{{\bf Z}}
\def\SL{{\rm SL}}
\def\PGL{{\rm PGL}}
\def\Aut{{\rm Aut}}
\def\odd{{\rm odd}}
\def\scirc{{\scriptstyle\circ}}
\def\emph#1{{\bf #1}}

\def\C{{\bf C}}
\def\F{{\bf F}}
\def\Z{{\bf Z}}
\def\N{{\bf N}}
\def\GL{{\rm GL}}
\def\End{{\rm End}}

\def\emph#1{{\it #1}}
\def\colon{:\,}
\def\SL{{\rm SL}}
\def\Ind{{\rm Ind}}
\def\SU{{\rm SU}}
\def\SO{{\rm SO}}

\def\boxotimes{\otimes}
\def\scirc{{\scriptstyle\circ}}
\def\G{G}
\def\isom{\cong}

\def\m#1#2#3#4{\left(\matrix{#1&#2\cr #3&#4\cr}\right)}
\def\F{{\bf F}}
\def\Z{{\bf Z}}
\def\SL{{\rm SL}}
\def\PGL{{\rm PGL}}
\def\Aut{{\rm Aut}}
\def\odd{{\rm odd}}
\def\scirc{{\scriptstyle\circ}}
\def\emph#1{{\bf #1}}

\centerline{\bf Density of the $\SO(3)$ TQFT representation of
mapping class groups}
\bigskip
\centerline{\byabs by}
\medskip
\noindent{\author\hfill Michael Larsen\footnote*{\sevenrm
Partially supported by NSF DMS 0100537 and DMS 0354772.} and Zhenghan Wang\footnote+{\sevenrm
Partially supported by NSF EIA 0130388, and DMS 0354772, and ARO.}\hfill}
\medskip
\centerline{\addr Department of Mathematics, Indiana University}
\centerline{\addr Bloomington, IN 47405, USA\footnote{}{AMS 2000
Classification 57R56}} \centerline{\email larsen@math.indiana.edu
and zhewang@indiana.edu}
\medskip

\def\C{{\Bbb C}}
\def\Z{{\Bbb Z}}

\def\N{{\Bbb N}}

\def\qed{\hfill\nobreak\rlap{$\sqcup$}$\sqcap$}

\def\F{{\Bbb F}}

\noindent {\bf \S 0. Introduction}
\medskip

A (2+1)-dimensional topological quantum field theory (TQFT) determines,
for each $g\ge 0$, a projective representation $(\rho_g,V_g)$ of the mapping
class group $M_g$ of a closed oriented surface of genus $g$.  This paper
is concerned with the $\SO(3)$ TQFT at an $r$th
root of unity, $r\ge 5$ prime.  Those TQFTs were first constructed mathematically
in [T].  The problem we consider is this: what is
the closure of $\rho_g(M_g)$?

For $g=1$ and the $\SU(2)$-theory, Kontsevich observed that the image is
finite.  A proof of finiteness both for the $\SU(2)$-theory and the $\SO(3)$-theory
may be found in [G]; see also [J], for an early calculation from
which the finiteness can be deduced.  For $g\ge 2$, the image
was shown to be infinite, and its closure therefore of
positive dimension [F].  In this paper, we identify the
representation for $g=1$ and show that the image is either
$\SL_2(\F_r)$ or $\PSL_2(\F_r)$ depending on whether $r$
is congruent to $1$ or $-1$ $({\rm mod}\ 4)$.  For $g\ge 2$, on the other hand,
we show that the image is as large as possible---that is, it is dense in
the group of projective unitary transformations on the representation
space of $M_g$.

We are working with the $\SO(3)$-theory because the statements are
cleaner than for the $\SU(2)$-theory.  Using the density result
here and the tensor decomposition formulas in Theorem 1.5 [BHMV],
the closed images
of the $\SU(2)$-theory can be identified.  The restriction to a prime $r$
is dictated by a result of Roberts [R] of which we make essential use: the
$\SU(2)$-theory representations of $M_g$ are irreducible for
$r$ prime.

The case $r=3$ is trivial since $\dim V_g = 1$ for all $g$.
The case $r=5$ was treated in joint work with Freedman [FLW2].
Our proof depended on the observation that a Dehn twist
acts with $r-1\over 2$ distinct eigenvalues.  For $r=5$, this
means that the $M_g$ representations satisfy the ``two-eigenvalue
property.''  The compact Lie groups $G$ admitting a representation with
respect to which a generating conjugacy class has only two eigenvalues
can be classified, and the few possibilities can be reduced to one
by an examination of the branching rules for the restriction $M_{g-1}\subset M_g$
and dimension computations (especially the Verlinde formula).
For $r\ge 7$, we no longer have the two-eigenvalue property, and the number
of possibilities for $G$ grows rapidly with $r$.  Mainly, therefore, we
depend on the branching rules.  A crucial point is to prove that
the representations are tensor-indecomposable, i.e., not equivalent
to a tensor product of representations of lower degrees; this is precisely why
the $\SO(3)$ case is simpler than that of $\SU(2)$.  The reason tensor-indecomposability
is so important is that, coupled with irreducibility, it implies that the identity
component of the closure of $\rho_g(M_g)$ is a simple group, and this greatly
shortens the list of possibilities.

The original motivation for this work was topological quantum computation
in the sense of [FKLW], [FLW1], and [FLW2].  As in [FLW2], there are also
applications to the distribution of values of $3$-manifold invariants.
As a simple example, we show that the set of the norms of the
Witten-Reshitikhin-Turaev $\SO(3)$ invariants of all connected 3-manifolds at
$A=ie^{2\pi i \over 4r}$ is dense in $[0,\infty)$ for primes $r\geq 5$.

\bigskip
 \noindent {\bf \S 1. The $\SO(3)$-TQFT}
\bigskip

There are several constructions of the $\SU(2)$ and the $\SO(3)$
TQFTs in the literature (e.g., [BHMV][FK][RT][T]).  The $\SU(2)$ TQFT
was first constructed mathematically in [RT], and the $\SO(3)$ TQFT in [T].
We will
follow Turaev's book [T], where the construction of a TQFT
is reduced to the construction of a modular tensor category.

Fixing a prime $r\geq 5$ and setting $A = ie^{2\pi i\over 4r} =
e^{2\pi i (r+1)\over 4r}$, note that $A$ is a primitive $2r$-th
root of unity when $r\equiv 1\; \hbox{mod} \; 4$, and a primitive $r$th
root of unity when $r\equiv -1\; \hbox{mod} \; 4$.  In [BHMV]
to construct TQFT using the
skein theory, the Kauffman variable $A$ is either a primitive $4r$th or a primitive $2r$th
root of unity.  When $r\equiv 1\; \hbox{mod} \; 4$ by the
$SO(3)$ TQFT we mean the TQFT denoted by $V_r$ in [BHMV] with the above choice
of $A$.
When $r\equiv -1\; \hbox{mod} \; 4$,  the same construction still gives rise
to a TQFT although $A$ is only
a primitive $r$th root of unity,
which is also denoted by $V_r$ here, but the decomposition formula in
Theorem 1.5 [BHMV] does not necessarily hold.  Consequently we have to
distinguish between the two cases $r\equiv 1 \; \hbox{mod} \; 4$
and $r\equiv -1 \; \hbox{mod} \; 4$.

The modular tensor categories associated to the TQFTs $V_r$ are described in
[T, Chapter XII].  In particular, [T, Theorem 9.2] discussed
the unitarity of the TQFTs.  For the above choice of $A$'s,  the ribbon categories
 of [T, Theorem 9.2] are not modular because the $S$-matrices as given
 in [T, Lemma 5.2] are
 singular. (The Kauffman variable $A$ is a primitive
$4r$th root of unity for only even $r$'s.)  But it can be shown that the even
subcategories (see [T, Section 7.5]) are indeed modular and unitary [T][FNWW].  The even
subcategories correspond to the restriction of
the representation categories to the odd dimensional (or integral spin)
\lq\lq halves"
in the quantum group setting [FK].

A modular tensor category consists of a large amount of data.
For our purpose here we will only specify the isomorphism classes of simple objects,
called {\it labels} of the associated TQFT, the $S$-matrix, and
the $T$ matrix.  More information is contained in Lemma 2.

We write the quantum integers
$[k]_A ={{A^{2k}-A^{-2k}}\over {A^2-A^{-2}}}$.  The label set of the $V_{r}$
theory is $L=\{0,2,4,\ldots,r-3\}$.  The quantum dimension of the
label $i$ is given by $d_i=[i+1]$, the subscript A in $[k]_A$ will be dropped from
now on,
and the global dimension of the $V_r$
modular tensor category is $D=\sqrt{\sum_{i \in L}[i+1]^2}={\sqrt{r} \over
2\sin{\pi \over r}}$.  The $\tilde{S}$-matrix $\tilde{S}=(\tilde{s}_{ij})$ can
be read off from Lemma 5.2 [T] as $\tilde{s}_{ij}=[(i+1)(j+1)]$.  The
$T=(t_{ij})$ matrix is diagonal with diagonal
entries the twists $\theta_{i}$, which are computed in [KL, Proposition 6, p. 43]
as  $\theta_i=A^{i(i+2)}$.

Let
$s=\m 0{-1}10$, $t=\m 1101$ be the generators of $\SL_2(\Z)$.
It is a deep fact that the $S, T$ matrices give rise to a projective
matrix representation
of $\SL_2(\Z)$ if we make the following assignments: $s\rightarrow
S={1\over D}\cdot \tilde{S}, t\rightarrow T^{-1}$, i.e., the
$SO(3)$ TQFT representation $\rho_{\SO(3)}$ for $\SL_2(\Z)$ is:
$$\rho_{\SO(3)}(s)={1\over D}([i+1][j+1])$$
and
$$\rho_{\SO(3)}(t)=( A^{-j(j+2)}\delta_{ij}).$$

The $T$-matrix corresponds geometrically to a Dehn twist, so
 the negative twists $\theta^{-1}_i$ are the eigenvalues for the image
of any Dehn twist on a non-separating simple closed curve.

\medskip
{\it Remark:} It seems to be generally believed that
the two theories $V_r$ and $V_{2r}$ constructed in [BHMV]
correspond to the $\SO(3)$ and the $\SU(2)$
Witten-Reshetikin-Turaev TQFTs. Actually the S-matrix of the $V_{2r}$
theory is not the same as that in the Witten-Reshetikin-Turaev $\SU(2)$ theory [Wi][RT]:
the $(i,j)$-th
entry differs by a sign $(-1)^{i+j}$.  But this discrepancy
disappears on restriction to the even subcategories; this is
the reason that the
Witten-Reshetikhin-Turaev TQFTs are always unitary, but the $V_{2r}$
theories are not unitary in general.  This subtle point is due to the
Frobenius-Schur indicators for self-dual representations, and will
be clarified in [FNWW].

\medskip

\nlemma\splitting{If $A$ is a primitive $2r$th root of unity
and $r \geq 3$ is odd, then there exists a TQFT $V_{2}'$
 and a natural isomorphism of theories such that
$$V_{2r}(\Sigma)\cong V_{2}'(\Sigma)\otimes V_{r}(\Sigma).$$
Moreover, the $\SO(3)$-theory representations of the
mapping class groups $M_g$ are irreducible for all primes $r\geq 5$.}

\bigskip

\proof The decomposition formula is Theorem 1.5 of [BHMV].

To prove the second part, first consider the case $r\equiv 1\; \hbox{mod} \; 4$.
Suppose the $\SO(3)$-theory representation of
the mapping class groups for a closed surface is reducible for a
prime $r$.  Then by the tensor decomposition formula above the
$\SU(2)$-theory representation would be reducible, too.  But this
contradicts the result of [R]. Therefore, the $\SO(3)$-theory
representations are also irreducible.  For $r\equiv -1\; \hbox{mod} \; 4$, the
same argument will work if a similar decomposition formula holds.
Without such a formula, the irreducibility of the $SO(3)$ representations
of $M_g$ can be deduced by following Roberts's argument [R].\qed

\nlemma\dimensions{Let $d_{r,g,h_1,\ldots,h_k}$ denote the
dimension of the vector space associated by the ${\rm SO}(3)$
theory at an $r$th root of unity to a compact
 oriented surface $\Sigma$ of genus $g$ with $k$
boundary components labeled $h_1,\ldots,h_k\in 2\Z$.
Then we have the following:
\smallskip
\item{1)}$d_{r,0,h}=\delta_{h 0}$, where $\delta$ is the Kronecker delta.
\item{2)}$d_{r,0,h_1,h_2} = \delta_{h_1 h_2}$.
\item{3)}$d_{r,0,h_1,h_2,h_3}$ is $1$ if and only if the $h_i$ satisfy the triangle
inequality (possibly with equality) and $h_1+h_2+h_3 \le 2(r-2)$; otherwise
it is $0$.
\item{4)}(Gluing formula) Suppose $\Sigma$ is cut along
a simple closed curve, and $\Sigma_{x}$ denotes the resulting surface
with the two new boundary
components both labeled by $x$.  Then
$$ V_r(\Sigma)\cong \oplus_{x\in L}V_{r}(\Sigma_{x}),$$ where
$L=\{0,2,\cdots, r-3\}$ is the label set of the $V_{r}$ theory.
\item{5)}$d_{r,1,h} = {r-1-h\over 2}$.
\item{6)}$d_{r,g,h_1,\ldots,h_k,0} = d_{r,g,h_1,\ldots,h_k}$.
\item{7)}(Verlinde formula) $d_{r,g} = \sum_{j=1}^{r-1\over 2}\alpha_j^{g-1}$, where
$\alpha_j = {r\csc^2 {2\pi j\over r}\over 4}$.
}

\proof Parts 1)--4) are the basic data of the $V_r$ theory (see also
[KL]).
Parts 5) and 6) are easy consequences of 1)--4).  The Verlinde formula
is derived in [BHMV]. \qed

\bigskip
\noindent {\bf \S 2. Tensor products and decompositions}
\medskip

Next we prove some technical results which enable us to establish
that certain representations are tensor indecomposible.  We say a
complex representation $V$ of a compact Lie group $G$ is \emph{isotypic}
if it is of the form $W^n = W\otimes\C^n$ for some irreducible representation
$W$ of $G$.   If $W$ is one-dimensional, we say $V$ is \emph{scalar}.
Two representations of $G$ are \emph{conjugate} if one is equivalent
to the composition of the other with an automorphism of $G$.

\medskip

\nlemma\decomp{Let $0\to G_1\to G_2\to G_3\to 0$ be a
short exact sequence of compact Lie groups and $\rho\colon G_2\to
\GL(V)$ an irreducible representation.  If
the restriction of $V$ to $G_1$ is isotypic, then there
exist central extensions $\tilde G_1$ and $\tilde G_2$ of $G_1$ and $G_2$
respectively and representations
$\sigma$ and $\tau$ of $\tilde G_2$ such that
\smallskip
\item{(1)} the extension $\tilde G_1$ is a normal subgroup of $\tilde G_2$,
\item{(2)} the quotient $\tilde G_2/\tilde G_1$ is isomorphic to $G_3$,
\item{(3)} the restriction of $\sigma$ to $\tilde G_1$ is irreducible,
\item{(4)} the restriction of $\tau$ to $\tilde G_1$ is scalar, and
\item{(5)} the tensor product $\sigma\otimes\tau$ is equivalent to
    the composition of $\rho$ with the central quotient map $\tilde G_2\to G_2$.}

\proof By hypothesis, $V|_{G_1}$ is isotypic and so
can be written as $W\otimes\C^k$.
The span of $\rho(G_1)$ is
$$\rho(G_1)\C=\End(W)\subset\End(W\otimes\C^k)=\End(V),$$
where $\End(W)$ maps to $\End(W\otimes\C^k)$
by $x\mapsto x\otimes {\rm Id}_k$.  The image $\rho(G_2)$ lies in the normalizer
$\End(W)\End(\C^k)$ of $\rho(G_1)\C$.  Thus, $\rho$ can be
regarded as a map
$$G_2\to(\End(W)\End(\C^k))^*=(\GL(W)\times\GL_k(\C))/\C^*.$$
Let
$$\tilde G_2 = G_2\times_{(\GL(W)\times\GL_k(\C))/\C^*} (\GL(W)\times\GL_k(\C)),$$
be a central extension of $G_2$, $\tilde G_1$ the pre-image of $G_1$ in
$\tilde G_2$ with respect to the central quotient map $\pi\colon \tilde G_2\to G_2$,
and $\tilde\rho$ the pullback $\tilde G_2\to  (\GL(W)\times\GL_k(\C))$
of $\rho$.  Let $\sigma$ and $\tau$ denote the compositions of $\tilde\rho$
with the projection maps $\GL(W)\times\GL_k(\C)\to\GL(W)$ and
$\GL(W)\times\GL_k(\C)\to\GL_k(\C)$ respectively.
The diagram
$$\matrix{
G_1&\leftarrow&\tilde G_1&\to&\tilde G_2\cr
\downarrow&&\downarrow&&\downarrow\rlap{$\scriptstyle\tilde\rho$}\cr
GL(W)= (\GL(W)\times\C^*)/\C^*&\leftarrow&\GL(W)\times\C^*&\to&
\GL(W)\times\GL_k(\C)\cr
}$$
shows that the restrictions of $\sigma$ and $\tau$ to $\tilde G_1$ are irreducible
and scalar respectively, and (5) is immediate.
\qed
\medskip
In the other direction, we have:
\medskip

\nlemma\product
{Let $G_2$ be a compact Lie group and $G_1$ a closed normal
subgroup.  Let $V$ and $W$ be irreducible representations of $G_2$
such that $V|_{G_1}$ is irreducible and $W|_{G_1}$ is scalar.
Then $V\otimes W$ is irreducible.}

\proof
As $V$ and $W$ are irreducible,
$$\dim\End_{G_2}(V)=\dim\bigl(V\otimes V^*\bigr)^{G_2}=1;
\ \dim\End_{G_2}(W)=\dim\bigl(W\otimes W^*\bigr)^{G_2}=1.$$
Let
$$V\otimes V^*=\bigoplus_i V_i,\ W\otimes W^*=\bigoplus_j W_j$$
denote decompositions into irreducible $G_2$-representations,
numbered so that $V_1$ and $W_1$ are trivial (and therefore the
other $V_i$ and $W_j$ are non-trivial). Thus,
$$\End_{G_2}(V\otimes W)=\bigl(\bigl(V\otimes V^*\bigr)
\otimes\bigl(W\otimes W^*\bigr)\bigr)^{G_2}
=\bigoplus_{i,j}(V_i\otimes W_j)^{G_2} =\bigoplus_{\{i,j\mid
V_i\isom W_j^*\}}\C.$$ Now, $W_j|_{G_1}$ is trivial for all $j$,
so $V_i\isom W_j^*$ implies that $V_i|_{G_1}$ is trivial.
However, $$1=\dim\End_{G_1}(V)=\dim\bigl(V\otimes
V^*\bigr)^{G_1}=\sum_i\dim V_i^{G_1},$$ so this is possible only
for $i=1$.  Then $W_j\isom V_1^*$ implies $j=1$, so
$$\dim\End_{G_2}(V\otimes W)=1.$$
\qed

\medskip
\nlemma\ltwo{For $r\ge 5$, the tensor product of any
two non-trivial irreducible representations of $\SL_2(\F_r)$ has
an irreducible factor of degree $>{r-1\over 2}$.}

\proof This can be deduced from the character table
[Sp]\ p. 160, whose notation we follow.    Let $\chi_1$ and
$\chi_2$ be non-trivial irreducible characters of $\SL_2(\F_r)$.
If $$\chi_1\chi_2=e\chi_\beta^++f\chi_\beta^-+g,$$ then
$$g=\cases{1&if $\chi_2=\bar\chi_1$,\cr 0&otherwise.}$$ If $d_b\in
T_1$ is a square, then $\chi_\beta^{\pm1}(d_b)=-1$.
If $\chi_2=\bar\chi_1$, comparing values at $d_b$,
$e+f\le 1$, which is absurd. Thus $g=0$, and thus
$|\chi_1(d_b)\chi_2(d_b)|=e+f$. As $\chi(1)\ge{r-1\over
2}|\chi(d_b)|$ for all non-trivial $\chi$, we get a contradiction.\qed

\nlemma\induce{Consider a short exact sequence of compact Lie groups
$$0\to G_1\to G_2\to \PSL_2(\F_r)\to 0,$$
where $r\ge 7$.  If
$H_2$ is the inverse image of a proper subgroup $H\subset \PSL_2(\F_r)$,
and $V$ is a representation of $H_2$, then some irreducible factor of $\Ind_{H_2}^{G_2} V$
has degree $> {r-1\over 2}$.}

\proof
Without loss of generality, we may assume that $V$ is irreducible.
Let $K = \ker G_1\to\GL(\Ind_{H_2}^{G_2} V)$.  Replacing $G_2$ and $H_2$
by $G_2/K$ and $H_2/K$ respectively if necessary,
we may assume the restriction of $\Ind_{H_2}^{G_2} V$ to $G_1$ is faithful.
Suppose that the restriction of some irreducible factor $W$ of $\Ind_{H_2}^{G_2} V$ to $G_1$
fails to be isotypic.    By Clifford's theorem, $W$
is a direct sum of mutually conjugate isotypic representations $W_i^k$ of $G_1$; the stabilizer of
$W_1^k$ can be regarded as a subgroup $\Delta$ of $\PSL_2(\F_r)$, and
$$\dim W = k \dim W_1 [\PSL_2(\F_r):\Delta].$$
By the well-known classification of subgroups of $\PSL_2(\F_r)$, each proper subgroup has index
$> {r-1\over 2}$, so $\dim W > {r-1\over 2}$.  We may therefore assume that the center of $G_1$ is in the center
of $G_2$.  By a theorem of Eilenberg and Mac Lane [EM], the obstruction to finding a
section of $G_2\to\PSL_2(\F_r)$ lies in $H^2(\PSL_2(\F_r), Z(G_1))$.  As $\PSL_2(\F_r)$ is perfect with universal central extension $\SL_2(\F_r)$, $G_2$ contains a subgroup $\Gamma$ isomorphic to $\PSL_2(\F_r)$
or $\SL_2(\F_r)$ which maps onto $\PSL_2(\F_r)$.  As
$$\Res_{G_2}^\Gamma \Ind_{H_2}^{G_2} V = \Ind_{H_2\cap\Gamma}^\Gamma\Res_{G_2}^\Gamma V,$$
we can reduce to the case that  $G_2$ is $\PSL_2(\F_r)$ or $\SL_2(\F_r)$.
If $G_2 = \PSL_2(\F_r)$ and $\tilde H_2$ is the preimage of $H_2\subset G_2$ in $\SL_2(\F_r)$,
then $\Ind_{H_2}^{G_2} V$, regarded as a representation of $\SL_2(\F_r)$ is the same as
$\Ind_{\tilde H_2}^{\SL_2(\F_r)} \tilde V$, where $\tilde V$ is $V$ regarded as a representation of $\tilde H_2$.
Therefore, without loss of generality, we may (and do) assume $G = \SL_2(\F_r)$.

If $H_1\subset H_2\subset G$,
$$\Ind_{H_2}^{G}\Ind_{H_1}^{H_2} V_1 = \Ind_{H_1}^{G} V_1,$$
so without loss of generality we may assume $H$ is a maximal proper subgroup.
Aside from the trivial representation, $\SL_2(\F_r)$ has two irreducible representations
of dimension $\le {r-1\over 2}$, with characters $\chi_\beta^\pm$.  As $V$
is irreducible, it is a subrepresentation of the regular representation of $H$,
and it follows that $\Ind_H^G V$ is a subrepresentation of the regular representation
of $\SL_2(\F_r)$.  In particular, $[G:H]\dim V$ reduces to 0 or 1 (mod ${r-1\over 2}$),
the former if $\dim V > 1$,
and
$$[G:H]\dim V\le 1^2 +( \dim\chi_\beta^+)^2+( \dim\chi_\beta^-)^2 = {r^2-2r+3\over 2}
< {|G|\over 2(r+1)}.$$
It follows that $|H| > 2(r+1)$.  By the classification of maximal subgroups of $\SL_2(\F_r)$,
this means that $H$ is a Borel subgroup $B$ of $G$ or
that $|H|\in\{24, 60, 120\}$.
The irreducible representations of $B$
all have degree $1$ or $r-1$.   If $\dim V=1$, the induced representation has degree
$r+1$, which does not satisfy the congruence condition.  If $\dim V = r-1$, the induced
representation has degree $r^2-1$,
which does not satisfy the inequality condition.
This leaves a short list of possible triples $(r,\dim V,|H|)$.   For $r > 11$, all can
be ruled out by the congruence condition or the inequality condition.
The only triples not ruled out are $(7,1,48)$ and $(11,1,120)$.
In each case, the degree of the induced representation is congruent to $1$ (mod ${r-1\over 2}$),
so $V$ must be trivial.   By [At]~pp.~3,~7, the induced representation in each case has an irreducible factor of degree $r-1$.\qed
\qed

\medskip
\nlemma\lquot{Consider a short exact sequence of compact
Lie groups
$$0\to G_1\to G_2\to \Gamma\to 0,$$
where $\Gamma$ is $\SL_2(\F_r)$ or $\PSL_2(\F_r)$, $r\ge 5$.
Suppose $V$ and $W$ are representations of $G_2$. If $V\otimes W$ has
\smallskip
\item{(1)}all its irreducible factors of  degree
$\le{r-1\over 2}$,
\item{(2)}at least one irreducible factor of degree $r-1\over 2$ which is $G_1$-scalar.
\item{(3)}exactly one irreducible factor of degree $1$,
\smallskip
then either $V$ or $W$ is one-dimensional.}

\proof Suppose first that $V$ and $W$ are irreducible and $V\otimes W$ satisfies
hypothesis (1).
By Clifford's theorem, we can write $V|_{G_1}$ and $W|_{G_1}$ as direct
sums of mutually conjugate isotypic representations $V_i^m$ and
$W_j^n$ respectively.
Let $H_2\supset G_1$ denote the subgroup of $G_2$ stabilizing both $V_1^m$ and $W_1^n$.
Then $\Ind_{H_2}^{G_2} V_1^m\otimes W_1^n$ is a $G_2$-subrepresentation of $V\otimes W$.
By \induce, it follows that $V|_{G_1}\isom V_1^m$ and $W|_{G_1}\isom W_1^n$.   (More generally,
if any tensor product of $G_2$-representations satisfies (1), all of the irreducible constituents of
all of the tensor factors are $G_1$-isotypic.)

By \decomp, replacing $G_1$ and $G_2$ by central extensions if necessary,
we can write $V$ as $V_\sigma\otimes V_\tau$
and $W$ as $W_\sigma\otimes W_\tau$ with $V_\sigma$ and
$W_\sigma$ $G_1$-irreducible and  $V_\tau$  and
$W_\tau$ $G_1$-scalar.   As explained above, any irreducible consistituent $U$ of
$V_\sigma\otimes W_\sigma$ is isotypic for $G_1$;
passing to central extensions of $G_1$ and $G_2$ if necessary,
we write $U$ as $U_\sigma\otimes U_\tau$, so $U_\sigma\otimes (U_\tau\otimes V_\tau\otimes W_\tau)$
is a $G_2$-subrepresentation of $V\otimes W$.
Every irreducible factor of $U_\tau\otimes V_\tau\otimes W_\tau$ is $G_1$-scalar.  By
\ltwo, unless at least two of $U_\tau$, $V_\tau$, and $W_\tau$ have dimension $1$, their tensor product
has an irreducible factor of degree $>{r-1\over 2}$ which is $G_1$-scalar.
If $\dim V_\tau = \dim W_\tau = 1$, then $V|_{G_1}$ and $W|_{G_1}$ are irreducible;
if $X$ is any one-dimensional representation of $G_1$, then $X^*\otimes V|_{G_1}$
is irreducible, so $\dim (X^*\otimes V\otimes W)^{G_1} \le 1$.

Thus, if $V\otimes W$ additionally satisfies hypothesis (2), either $\dim V_\tau > 1$ or $\dim W_\tau > 1$.
Without loss of generality, we assume the latter is true.   Unless $\dim U_\tau = \dim V_{\tau} = 1$,
$V\otimes W$ contains a factor which is the tensor product of the $G_1$-irreducible $U_\sigma$ with
a $G_1$-scalar irreducible of dimension $>{r-1\over 2}$.  This is impossible by \product.
The situation is therefore that $V_\tau\otimes W_\tau$ is $G_1$-scalar, $G_2$-irreducible and of
dimension ${r-1\over 2}$, and
and every factor of $V_\sigma\otimes W_\sigma$ is $G_1$-irreducible.  Applying (2) and \product\ again,
we see further that every $G_2$-irreducible factor of $V_\sigma\otimes W_\sigma$ has dimension 1, so
$W_\sigma^*\otimes W_\sigma$ decomposes entirely into $1$-dimensional pieces over $G_2$.
If $V'$ is any other irreducible representation such that $V'\otimes W$ satisfies (1), then $V'_\sigma$ is
a twist of $W_\sigma^*$ by a character so $V'_\sigma\otimes {V'_\sigma}^*$ also decomposes as a sum of
$(\dim V'_\sigma)^2$ characters over $G_2$.

Now we return to the original problem.  If $V\otimes W$ satisfies all three hypotheses, there
exists irreducible factors $V_0$ and $W_0$ of
$V$ and $W$ respectively satisfying hypotheses (1) and (2).
Without loss of generality, we assume $\dim V_{0\tau}=1$.
Then every irreducible
factor of $V$ is a twist of $W^*_{0\sigma}$, so they all have the same dimension $d$.
By hypothesis (3), there exist irreducible factors $V_i$ and $W_j$
of $V$ and $W$ respectively such that $V_i\otimes W_j$ has a one-dimensional
$G_2$-submodule.   Thus $W_j^*$ must be a twist
of $V_i$, and $(V\otimes W_j)|_{G_1}$
decomposes entirely into $d^2$ $1$-dimensional pieces.  By (3), $d=1$.

\qed

\bigskip
\noindent {\bf \S 3. Case g=1}

\bigskip

\noindent {\bf Theorem 1.} The projective representation of
$M_1 = {\rm Map}(S^1\times S^1) = \SL_2(\Z)$
given by the $\SO(3)$-theory at $A=ie^{2\pi i\over 4r}$ is the same
as the projective representation obtained by composing the (mod
$r$) reduction map $\SL_2(\Z)\to \SL_2(\F_r)$ with the odd factor
of the Weil representation of $\SL_2(\F_r)$.
\bigskip

We want explicit matrices for one of the two ${r-1\over 2}$-dimensional
irreducible representations of $\SL_2(\F_r)$.  To find them,
we briefly recall the theory of Weil representations over finite
fields [Ge].  Let $H_r$ denote the
Heisenberg group of order $r^3$. We regard $H_r$ as a central
extension
$$0\to\F_r\to H_r\to \F_r^2\to 0.$$
The extension class defines a symplectic form (in this case, an
area form) on the quotient.  Any automorphism of $H_r$ stabilizes
the center and acts on $\F_r^2$, respecting this symplectic form.
Regarding  $\SL_2(\F_r)$ as the group of symplectic linear
transformations of $\F_r^2$, we claim that its action lifts
to $H_r$. To make this explicit, let $x$ and $y$ be elements of
$H_r$ whose images in $\F_r^2$ form a unimodular basis.  Let $z$ be
the generator of the center defined by $z^2 = yxy^{-1}x^{-1}$.
Finally, for $M = {\m abcd}\in \SL_2(\F_r)$, let
$$f_M(x) = z^{ac}x^a y^c,\ f_M(y) = z^{bd}x^by^d,\ f_M(z)=z.$$
We easily check that this defines an action of $\SL_2(\F_r)$ on $H_r$.

Let $e(k)$ denote $e^{2\pi i k\over r}$. Let $(\rho,V)$ denote the
Stone-von Neumann representation (i.e., the unique irreducible
representation of $H_r$ with central character $z^k\mapsto e(k).$)
We fix a basis $e_0,\ldots,e_{r-1}$ for $V$ so that $e_k$ is the
$e(2k)$-eigenspace of $x$. In this basis, we write
$$\rho(x)=\left(\matrix{ 1&0&0&\cdots&0\cr 0&e(2)&0&\cdots&0\cr
0&0&e(4)&\cdots&0\cr \vdots&\vdots&\vdots&\ddots&\vdots\cr
0&0&0&\cdots&e(-2)\cr }\right),\quad \rho(y)=\left(\matrix{
0&1&0&\cdots&0\cr 0&0&1&\cdots&0\cr
\vdots&\vdots&\vdots&\ddots&\vdots\cr 0&0&0&\cdots&1\cr
1&0&0&\cdots&0\cr }\right),\quad \rho(z) = e(1)\,{\rm Id}. $$
For any $\alpha\in \SL_2(\F_r)\subset\Aut(H_r)$, $\rho\scirc\alpha$ is
equivalent to $\rho$. There exists $R_\alpha$, therefore, unique
up to scalar multiples, such that
$$R_\alpha \rho(h) R_\alpha^{-1} = \rho(\alpha(h))\eqno{(*)}$$
for all $h\in H_r$.  If $\bar R_\alpha$ denotes the class of
$R_\alpha$ in $\PGL(V)$, we conclude that $\alpha\mapsto \bar
R_\alpha$ is a projective representation. When $r\ge 5$, $\SL_2(\F_r)$ is
perfect and centrally closed, so there is a unique lifting to an
$r$-dimensional linear representation, which we call the
\emph{Weil representation} of $\SL_2(\F_r)$.
Explicitly we may choose (up to scalar multiplication)
$$R_S = (e(2ij))_{0\le i,j < r},\ R_T
= (\delta_{ij}e(-i^2))_{0\le i,j < r}.$$
We can verify ($\ast$) by checking it for $h=x$ and $h=y$.

Let $E=\{0,2,4,\ldots,r-3\}$ and set
$$f_i = e_{r-1-i\over 2} - e_{r+1+i\over 2},\ i\in E.$$
It is easy to see that the span of $f_0,f_2,f_4,\ldots,f_{r-3}$
forms an invariant subspace $V^{\odd}$ of both $R_S$ and $R_T$. In
terms of this basis, $R_S$ is represented by the matrix
$$R_S^{\odd} = \left( e\left(2\left({r-1-i\over
2}\right)\left({r-1-j\over 2}\right)\right) -
e\left(-2\left({r-1-i\over 2}\right)\left({r-1-j\over
2}\right)\right) \right)_{i,j\in E},$$
and $R_T$ is represented by
$$R_T^{\odd} = \left(\delta_{ij}e\left(-\left({r-1-j\over
2}\right)^2\right)\right)_{i,j\in E}.$$
As
$$\ A^2 = e\Bigl({r+1\over 2}\Bigr),$$
we obtain
$$R_S^{\odd} =
\left(A^{2(i+1)(j+2)}-A^{-2(i+1)(j+1)}\right)_{i,j\in E} =
(A^2-A^{-2})\bigl([(i+1)(j+1)]_A\bigr)_{i,j\in E},$$
and
$$R_T^{\odd}= e^{-{2\pi i(r-1)^2\over 4r}} \left(\delta_{ij}
A^{-j(j+2)} \right)_{i,j\in E}.$$

Thus the composition $\SL_2(\Z)\to\SL_2(\F_r)\to\PGL(V^\odd)$
is equivalent to $\rho_{\SO(3)}$.
\qed

\bigskip
\noindent {\bf \S 4. Case g=2}
\medskip

In what follows, we write $X_r$ for $$X_r=\left\{e^{2\pi i
n^2\over r}\Bigm| 0<n<{r\over 2}\right\}.$$
\smallskip

\medskip
\nlemma\cox{Let $G$ be a simple compact Lie group with
Coxeter number $h$, $r\ge 5$ a prime, and $\rho:\,G\to\GL(V)$ an
irreducible representation of $G$ such that $r$ divides $\dim V$.
If there exists $g\in G$ such that the spectrum of $\rho(g)$ is
$X_r$, then $r\le 2h-5$.}

\proof
Let $T$ be a maximal torus containing $g$ and $(\,,\,)$ denote the
Cartan-Killing form on the character space $X^*(T)\otimes{\bf R}$.  Let
$$\langle\beta,\alpha\rangle={2(\beta,\alpha)\over(\alpha,\alpha)}, $$
and fix a Weyl chamber.
If $V$ has highest weight $\lambda$ and
$\rho$ is the half sum of
positive roots, the Weyl dimension formula ([Bo]~VIII,~\S9,~Th. 2) asserts
$$\dim V=\prod_{\alpha>0}{\langle\lambda,
\alpha\rangle+\langle\rho,\alpha\rangle\over\langle\rho,\alpha\rangle},$$
where the product is taken over all positive roots.  Let $\beta$ denote
the highest root.  Since $r$ divides $\dim V$, and
$\langle\mu,\alpha\rangle\in\Z$ for all weights $\mu$ and roots $\alpha$,
$$\langle\lambda,\beta\rangle+\langle\rho,\beta\rangle\ge r, $$
or, by [Bo]~VI,~\S1,~Prop.~29(c) and [Bo]~VI,~\S1,~Prop. 31,
$$\langle\lambda,\beta\rangle+h-1\ge r.$$
The string of weights $\lambda,\lambda-\beta,\ldots$ has length
$1+\langle\lambda,\beta\rangle$ ([Bo]~VIII,~\S7,~Prop.~3(i)). For
any $w$ in the Weyl group, the string
$w(\lambda),w(\lambda)-w(\beta),\ldots$ has the same length.  The
Weyl orbit of $\beta$ consists of all long roots
([Bo]~VI,~\S1,~Prop.~11), so the lattice it generates contains the
root lattice if $G$ is of type A, D, or E; twice the root lattice
if $G$ is of type B, C, or F; and three times the root lattice if
$G$ is of type G ([Bo]~Planches). Since the difference between
weights in an irreducible representation belongs to the root
lattice, and since the eigenvalues of $\rho(g)$ are $r$th roots of
unity, not all equal, we conclude that if $r\ge 5$, $w(\beta)(g)$
is a primitive $r$th root of unity for some $w\in W$.  A
non-trivial geometric progression in $X_r$ has length $\le
{r-1\over 2}$, so
$$ 1+\langle\lambda, \beta\rangle\le{r-1\over 2}.$$
Thus, $r\le 2h-5$.\qed

\medskip
\nlemma\class{Under the hypotheses of \cox, if $\dim
V={r^3-r\over 24}$, then $G$ is a classical group.  If $V$ is not
self-dual, then $G$ is of type $A_n$.}

\proof The Coxeter numbers of $G_2$, $F_4$, $E_6$,
$E_7$, $E_8$ are $6$, $12$, $12$, $18$, $30$, respectively
([Bo]~Planches).   Examination of all primes $\le 55$ [MP] reveals
that the only case in which $r^3-r\over 24$ is the dimension of a
representation  of a suitable exceptional group is $r=7$.   The
group is $G_2$, and $V$ is the adjoint representation.  This case
is excluded, however, as the longest string of short roots is
$4>{r-1\over 2}$.  Thus $\beta(g)=1$ for all short roots $\beta$.
As the short roots generate the root lattice, $g$ has
all eigenvalues equal, contrary to assumption.

If $V$ is not self-dual, $G$ cannot be of type $B_n$ or $C_n$ and
can only be of type $D_n$ if $n$ is odd and the highest weight
$\lambda=a_1\lambda_1+\cdots+a_n\lambda_n$ satisfies
$\sup(a_{n-1}, a_n)>0$ [MP].  As the Weyl dimension formula is
monotonic in each $a_i$, $\dim V\ge 2^{n-1}$ while $r\le
2h-5=4n-9$.  Given the dimension of $V$, the only possibilities
for $(n,r)$ are $(5,11)$, $(7,13)$, $(7,17)$, $(7,19)$, $(9,19)$,
$(9,23)$, $(11,31)$.  For $r\le 31$, the longest geometric
progression in $X_r$ has length $\le 4$, so $h\ge r-2$ or $r\le
2n$.  This leaves only the case $(7,13)$, and by [MP], the only
irreducible $91$-dimensional representation of $D_7$ is self-dual.\qed

\medskip
\nlemma\lodim{For $n\ge 11$,
if $G=\SU(n)$ and $V$ is an irreducible representation of $G$ with highest
weight $\lambda$ and dimension $<2{n\choose 3}$, then $\lambda$ is one of
the following:
\medskip
\centerline{\vbox{\offinterlineskip
\halign{&\vrule#&
  \strut\quad\hfil#\quad\cr
\multispan5\hrulefill\cr
height2pt&\omit&&\omit&\cr
&Highest weight\hfil&&Dimension&\cr
height2pt&\omit&&\omit&\cr
\multispan5\hrulefill\cr
height2pt&\omit&&\omit&\cr
&$0$&&$1$&\cr
&$\lambda_1,\,\lambda_{n-1}$&&$n$&\cr
&$\lambda_2,\,\lambda_{n-2}$&&$n\choose 2$&\cr
height2pt&\omit&&\omit&\cr
&$2\lambda_1,\,2\lambda_{n-1}$&&$n+1\choose 2$&\cr
height2pt&\omit&&\omit&\cr
&$\lambda_1+\lambda_{n-1}$&&$n^2-1$&\cr
&$\lambda_3,\,\lambda_{n-3}$&&$n\choose 3$&\cr
height2pt&\omit&&\omit&\cr
&$3\lambda_1,\,3\lambda_{n-1}$&&$n+2\choose 3$&\cr
height2pt&\omit&&\omit&\cr
\multispan5\hrulefill\cr}}}
}

\proof
By the monotonicity of the Weyl dimension formula, it is enough to check
that the dimension is always greater than $2{n\choose 3}$ in the following
cases (and their duals): $\lambda_k$, $4\le k\le n-4$; $4\lambda_1$;
$2\lambda_2$; $2\lambda_3$; $\lambda_1+\lambda_2$; $\lambda_1+\lambda_3$;
$\lambda_2+\lambda_3$; $\lambda_1+\lambda_{n-2}$;
$\lambda_1+\lambda_{n-3}$;
$\lambda_2+\lambda_{n-2}$; $\lambda_2+\lambda_{n-3}$;
$\lambda_3+\lambda_{n-3}$; $2\lambda_1+\lambda_{n-1}$.
\qed
\medskip

\nlemma\dense{Let $r\ge 7$ be a prime, $G$ a simple
compact Lie group, and $\rho:\,G\to\GL(V)$ an irreducible
representation of $G$ such that
\smallskip
\item{1)} $V$ is not self dual,
\item{2)} $\dim V={r^3-r\over 24}$,
\item{3)} There exists $g\in G$ such that the spectrum of $\rho(g)$ is
$X_r$.
\smallskip Then either $G=\SU(\dim V)$, and $V$ is the standard
representation or its dual; or $r=13$,  $G$ is a central quotient
of $\SU(14)$, and $V$ is the exterior square representation
or its dual.}

\proof By \cox, $r\le 2h-5$.  By \class, $G$ is
of type $A_n$, so $n\ge {r+3\over 2}$, so
$$2{n\choose 3}\ge{(r+3)(r+1)(r-1)\over 24}>{r^3-r\over 24}.$$
If $n+1 < 11$, $r < 17$, so there are three cases: $r = 7$, $n\ge 5$, and $\dim V = 14$;
$r = 11$, $n\ge 7$, and $\dim V = 55$; and $r = 13$, $n\ge 8$, and $\dim V = 91$.
For $n\le 10$, we see there is just one possibility for an irreducible representation
of $\SU(n+1)$ of the given dimension: $r = 11$, $n = 10$, and $V$ is the exterior square
of the standard representation of $\SU(11)$ or its dual.  If $n+1\ge 11$,
by \lodim, either $V$ or $V^*$ has highest weight in the set
$$\{\lambda_1,2\lambda_1,3\lambda_1, \lambda_2,\lambda_3,
\lambda_1+\lambda_{n-1}\}$$
As $V$ is not self-dual, we can exclude $\lambda_1+\lambda_{n-1}$.
For $r\ge 7$,
$${(r+1)r(r-1)\over 24}<{(r-1)(r-2)(r-3)\over 8}, $$
so ${m\choose 3}\le {r^3-r\over 24}$ only when $m<r-1$, in which
case equality is ruled out since $r$ does not divide ${m\choose
3}$. Applying this when $m=n\pm 1$, we exclude the cases
$\lambda_3$ and $3\lambda_1$. For $\lambda_2$ and $2\lambda_1$, we
seek solutions of
$${m\choose 2} = {r^3-r\over 24}$$
in integers $m$. For any such solution, $r\vert m$ or $r\mid m-1$,
so
$$12a(ar\pm 1) = r^2-1,\ a\in\N.$$
The discriminant of the quadratic equation for $r$ is $(12a^2)^2 +
4\pm 48a$, and
$$(12a^2-1)^2 < (12a^2)^2 + 4 - 48a< (12a^2)^2 < (12a^2)^2 + 4 + 48a <
(12a^2+1)^2$$
for $a\ge 3$.  For $a=2$, the discriminant is not square for
either choice of sign.  For $a=1$, we get the two solutions
$(r,m)=(11,11)$ and $(r,m)=(13,14)$.  An exhaustive analysis of
sets $S$ of $r$th roots of unity, whose symmetric or exterior
squares give $X_{11}$ or $X_{13}$ reveals exactly two
possibilities: the set
$$\{1,\zeta_{13},\zeta_{13}^3,\zeta_{13}^9\}\eqno{(1)}$$ and its
complex conjugate have exterior square $X_{13}$.\qed

Let $\rho_g:\,M_g\to \PGL(V_g)$ denote the projective unitary
representation given by the $\SO(3)$-theory.
Let $\G_g$ denote the
closure of the image.  It is a subgroup
of $\PSU(\dim V_g)$ and therefore a compact Lie group.
We will often regard $V_g$ as a projective
representation of $\G_g$.

\medskip
\noindent{\bf Theorem 2.} For $g=2$, the projective representation
$\rho$ associated to the $\SO(3)$ theory at $A=ie^{2\pi i\over 4r}$
for $r\geq 5$ has dense image.
\smallskip

The proof will be carried out in several steps.  By [FLW2],
we may assume $r\ge 7$.

\noindent{\bf Step 1.} The Lie group $\G_2$ is infinite.

\proof Consider the decomposition of the
representation space arising from a curve separating a genus $2$
surface into two genus $1$ surfaces with boundary.  The components
are indexed by labels $0,2,\ldots,r-3$, and they are
projective representation spaces of $M_1\times M_1$. The representation
associated with label $2l$ is of the form $W_{2l}\boxotimes
W_{2l}$, where each tensor factor has dimension $r-1-2l\over2$.
For label $r-5$, it has dimension 2.  Thus we have a
two-dimensional projective unitary representation of
$M_1=\SL_2(\Z)$; the ratio of eigenvalues for a Dehn twist is a
primitive $r$th root of unity. By the classification of finite
subgroups of $\SO(3)$, this implies the image is infinite, and it
follows that the same is true for $\G_2$.
\smallskip

\noindent{\bf Step 2.} The projective representation $V_2$ is not self-dual.

\proof Equivalently, for any central extension $\tilde M_2$ for which
one can lift $V_2$ to a linear representation (also denoted $V_2$),
the contragredient representation $V_2^*$ is not obtained by tensoring $V_2$
by a central character.
We compute the multiplicities of the
eigenvalues of a lift to $\tilde M_2$ of a Dehn twist. These are just the dimensions
$d_{r,1,2l,2l}$ of a doubly-punctured torus with both labels equal
to $2l$, and are therefore given by
$$d_{r,1,2l,2l} = {(2l+1)(r-2l-1)\over2}.\eqno{(3)}$$
No two of these multiplicities
coincide as $l$ ranges over integers $\le{r-3\over2}$, so $V_2$
cannot be self-dual.

\noindent{\bf Step 3.} Let $\tilde G_2$ denote any central extension
of $G_2$ for which $V_2$ lifts to a linear representation (which we also
denote $(\rho_2,V_2)$).  Let $\tilde G_2^\circ$ denote the identity
component of $\tilde G_2$.  Then the restriction of $V_2$ to
$\tilde G_2^\circ$  is isotypic.

\proof As $M_2$ is generated by Dehn twists, some
lift $\tilde t\in\tilde G_2$ of a Dehn twist $t$
would otherwise permute the isotypic components
non-trivially.  Thus, the eigenvalues of $\rho_2(\tilde t)$ (which are
defined up to multiplication by a common scalar) would contain a
coset of a non-trivial group of roots of unity [FLW2]~Lemma~1.2.
This is impossible, since up to scalars, the spectrum of a Dehn
twist is $X_r$.
\smallskip

\noindent{\bf Step 4.} For any central extension $\tilde G_2$ of $G_2$ as above and any
normal subgroup $\tilde G'_2$ such that every homomorphism
$\SL_2(\F_r)\to\tilde G_2/\tilde G'_2$ is trivial, $V_2$
is tensor indecomposable as a $\tilde G'_2$-representation.

\proof
The restriction of $V_2$ to
$$\tilde M_1=\widetilde{\{1\}\times M_1}\subset \widetilde{M_1\times M_1}\subset \tilde M_2$$
decomposes as a sum of terms of the form $W_{2l}$.  Let $\tilde H_1$ denote the
closure of $\tilde M_1$ in $\GL(V_2)$
and $\tilde H'_1$ the intersection of $\tilde H_1$ with $\tilde G'_2$.
The occurrence of $W_0$ as a
factor guarantees that $\SL_2(\F_r)$ or $\PSL_2(\F_r)$ is a quotient of $\tilde H_1$ for
which the $W_0$-factor in question is an irreducible representation of
degree $r-1\over 2$; the condition on $\tilde G'_2$ guarantees that
$\tilde H'_1$ maps onto $\SL_2(\F_r)$, so $V_2|_{\tilde H'_1}$
has an irreducible factor which is the composition of this quotient map
and a degree $r-1\over2$ representation of $\SL_2(\F_r)$.
Tensor indecomposability now follows from \lquot.
\smallskip

\noindent{\bf Step 5.} The restriction of $V_2$ to $\G_2^\circ$ is irreducible.

\proof
Otherwise, by \decomp, there exists a central extension $\tilde G_2$ of $G_2$
so that regarding $V_2$ as a $\tilde G_2$-representation, it has a tensor-decomposition.
\smallskip

\noindent{\bf Step 6.} The identity component $G_2^\circ$ is a simple compact Lie group.

\proof

The identity component  $K= G_2^\circ$ is a connected compact Lie group.
As $V_2$ is an irreducible projective representation, the center
of $K$ is trivial.  Therefore, it is
a product of compact simple Lie groups $K_1\times\cdots \times K_s$,
and $V_2$ is a tensor product of unitary projective representations $X_1,\ldots X_s$
of the $K_i$.   In other words,
$$K_1\times \cdots \times K_s\hookrightarrow \PSU(X_1)\times\cdots\times \PSU(X_s)
\hookrightarrow \PSU(V_2),
$$
where the first inclusion is the product of inclusion $K_i\hookrightarrow \PSU(X_i)$
and the second is the tensor product map.  As $2^r>{r^3-r\over24}$ for all $r$,
$s < r$.  Consider the composition $\pi\colon G_2\to\Aut(K)\to\Out(K)$.
The outer automorphism group is contained in a product of groups of the form
$\Out(K_i)^{s_i}\rtimes S_{s_i}$, where $s_i\le s$.
The largest proper subgroup of $\SL_2(\F_r)$ is the Borel subgroup with
index $r+1 > s_i$, so any homomorphism $\SL_2(\F_r)\to S_{s_i}$ is trivial.
Therefore, any homomorphism from $\SL_2(\F_r)$ to $\Out(K)$ lands in
a product of solvable groups, and since $\SL_2(\F_r)$ is perfect, that means any
such homomorphism is trivial.  By Step 4, $V_2$ is tensor indecomposable
as a representation of any central
extension of the subgroup $G'_2 = \ker\pi$.  However, $g$ acts by inner
automorphisms on $K$ for $g\in G'_2$, and since $\rho_2$ is irreducible,
this means $G'_2\subset K\subset \prod_{i=1}^s\PSU(X_i)$,
so $s=1$.

\smallskip

\noindent{\bf Step 7.} The theorem holds if $r\neq 13$.
\proof
Applying \dense\ to the universal cover of $G_2^\circ$,
$\G_2^\circ$ is all of $\PSU(\dim V_2)$, so
the same is true for $\G_2$.

\smallskip

\noindent{\bf Step 8.} The theorem holds if $r = 13$.
\proof

For $r=13$, we must consider the
possibility that the universal covering group of $\G_2^\circ$ is
$\SU(14)$, $V_2$ is its alternating square, and a Dehn twist has
exactly four different eigenvalues $\lambda_i$ in $\SU(14)$ given
up to a common scalar multiple by (1) or its complex conjugate. By
(3), the eigenvalues of a Dehn twist have multiplicities $6, 15,
20, 21, 18, 11$ in $\PSU(\dim V_2)$, and each one arises uniquely
as a product of distinct eigenvalues $\lambda_i$. Therefore, some
$\lambda_i$ must have multiplicity $11$, but this is impossible
since only one of the eigenvalues in $V_2$ has multiplicity
divisible by 11. Therefore, $\G_2^\circ=\PSU(\dim V_2)$ also for
$r=13$.

\smallskip
\qed

\bigskip
\noindent {\bf \S 5. Case g$\geq $3}
\medskip

\medskip
\noindent{\bf Theorem 3.} For all $r\ge 5$ and all $g\ge 2$,
$\rho_g(M_g)$ is dense in $\PSU(\dim V_g)$.
\smallskip
For $r=5$, this is already known [FLW2]~Theorem~6.2.
We therefore assume from now on
that $r\ge7$. We begin with a dimension estimate.

\medskip
\nlemma\goslow{For $r\ge7$ and $g\ge2$,
$$\dim V_{\Sigma_{g+1}}<{\dim V_{\Sigma_g}\choose 2}\eqno{(2)}$$
except  when $r=7$ and $g=2$.}

\proof As $\alpha_k>1$ for all $k$, if $g\ge 3$, $$
\dim V_{\Sigma_{g+1}}=\sum_k\alpha_k^g\le\sum_k
\alpha_k^{3(g-1)/2} <\dim^{3/2} V_{\Sigma_g}<{\dim
V_{\Sigma_g}\choose 2}.$$ For $g=2$, we compute $${\dim
V_{\Sigma_2}\choose 2}-\dim V_{\Sigma_3} =
{(r+5)(r+3)(r+1)r(r-1)(r-8)\over 5760},$$ and this quantity is
obviously positive when $r>7$.\qed

\smallskip

\noindent{\it Proof of Theorem 3.}

The proof is very similar to that of Theorem 2.

\noindent{\bf Step 1.} The Lie group $\G_g$ is infinite.

\proof
Consider the decomposition of the
representation space arising from a curve separating a genus $g$
surface into two pieces, one of genus 1 and one of genus $g-1$.
Restricting to $M_{g-1}\times M_1$, we obtain a decomposition
$$V_g=\bigoplus_{l=0}^{r-3\over2}X_{g-1,2l}\boxotimes W_{2l},$$
where $X_{g-1,2l}$ denotes the projective representation
space of $M_{g-1}$
associated to a surface of genus $g-1$ with a single boundary
component labeled $2l$.
Now we proceed as in Step 1 of Theorem 2.
\smallskip

\noindent{\bf Step 2.}
The projective representation $V_g$ is not self-dual.

\proof
If $V_g$ is self-dual, its restriction to
$M_{g-1}\times M_1$ decomposes into self-dual projective representations and
mutually dual pairs of projective representations.  We use induction on $g$,
the base case being Step 2 of Theorem 2.  By the induction
hypothesis, $X_{g-1,0} \boxotimes W_0 = V_{g-1}\boxotimes W_0$ is not self-dual.
Neither can it be dual to any other factor since
$W_0=V_1$ is irreducible and the other representations $W_{2l}$ have lower
dimension.
\smallskip

\noindent{\bf Step 3.}
Let $\tilde G_g$ denote any central extension
of $G_g$ for which $V_g$ lifts to a linear representation.
Let $\tilde G_g^\circ$ denote the identity
component of $\tilde G_g$.  Then the restriction of $V_g$ to
$\tilde G_g^\circ$  is isotypic.

\proof
Identical to the proof of Step 3 of Theorem 2.
\smallskip

\noindent{\bf Step 4.} For any central extension $\tilde G_g$
of $G_g$ as above and any normal subgroup, $V_g$
is tensor indecomposable as a $\tilde G_g$-representation.

\proof
The restriction of $V_g$ to
$$\tilde M_1=\widetilde{\{1\}\times M_1}\subset \widetilde{M_{g-1}\times M_1}
\subset \tilde M_g$$
decomposes as a sum of terms of the form $W_{2l}$.
Now we proceed as in Step 4 of Theorem 2.
\smallskip

\indent{\bf Step 5.} The restriction of $V_g$ to $G_g^\circ$ is irreducible.

\proof
Identical to the proof of Step 5 of Theorem 2.
\smallskip

\noindent{\bf Step 6.} The identity component $G_g^\circ$ is a simple compact Lie group.

\proof
Let $K=G_g^\circ$.  As in Step 6 of Theorem 2, $K$ is a product of compact
simple Lie groups $K_i$, and $V_g$ is the tensor product of
unitary projective representations $X_i$ of the $K_i$.  Thus,
$$K_1\times \cdots \times K_s\hookrightarrow \PSU(X_1)\times\cdots\times \PSU(X_s)
\hookrightarrow \PSU(V_g).
$$
The conjugation action of $G_2$ on $K$ must act transitively on the factors, since
any decomposition into orbits gives a tensor decomposition of $V$.  Therefore,
the $K_i$ are mutually isomorphic, and their representations $X_i$ are equivalent
up to composition with automorphisms of $X_i$; in particular, their degrees
are all the same.  Now,
the closure of $M_{g-1}\subset M_{g-1}\times M_1$ in $K$ maps onto $G_{g-1}$ since
$V_{g-1} = X_{g-1,0}$ is a summand of the restriction of $V_g$ to $M_{g-1}$.
Thus some factor $K_i$ maps onto $G_{g-1}$.  By the induction hypothesis,
$$\dim K_i\ge \dim G_{g-1} = \dim V_{g-1}^2-1,$$
so any non-trivial representation of $K_i$ has degree $\ge\dim V_{g-1}$.  Therefore,
$\dim V_g \ge \dim V_{g-1}^s$.  The inequality (2) then implies $s=1$ except possibly
when $g=3$ and $r=7$, in which case, $\dim V_g = 98 < 14^2 = \dim V_{g-1}^2$, so again $s=1$.

\smallskip

\noindent{\bf Step 7.} For all $g\ge 3$ and $r\ge 7$, $G_g = \PSU(\dim V_g)$.

\proof
We use induction, the base case being Theorem 2.
By the induction hypothesis and (2),
$${\rm rank}\,G_g^\circ\ge{\rm rank}\,G_{g-1}\ge\dim V_{g-1}-1>
\sqrt{2\dim V_g+1/4}-1/2.\eqno{(4)}$$
By Step 2, $V_g$ is not
self-dual, so by [MP], $G_g$ is of type $A_n$, $D_n$ ($n$ odd),
or $E_6$.  The case $E_6$ is ruled out since $\dim V_{g-1}>7$ in
all cases $g\ge3$.  The minimal dimension for a
representation of $D_n$ which is not self-dual is $2^{n-1}>{n+1\choose 2}$ for $n>4$.
This leaves the case $A_n$, where the inequality $\dim
V_g<{n+1\choose 2}$ implies $\dim V_g<2{n+1\choose3}$ for $n>2$.
\lodim\ gives the list of possibilities.  For $n>7$, the possible
highest weights are $\lambda_1$, $\lambda_n$, $\lambda_2$,
$\lambda_{n-1}$, $2\lambda_1$, $2\lambda_n$, and
$\lambda_1+\lambda_n$, the last being ruled out as $V_g$ is not
self-dual.  Up to duality, then, $V_g$ is either the standard
representation, its exterior square, or its symmetric square, and
both are ruled out by (4).

\bigskip
\noindent {\bf \S 6. An application}
\medskip

Besides determining the representations of the mapping class groups, a
TQFT also determines invariants of oriented closed 3-manifolds [RT].
To describe the $\SO(3)$ invariant of 3-manifolds, we introduce
the following notations:  let $d_i$ and $\theta_i$ be the
quantum dimension and the twist of the label $i$, and $D$ the global
dimension of the modular tensor category defined
in Section 1; then we define $p_{\pm}=\sum_{i\in L}\theta_i^{\pm
1} d_i^2$, and $\omega_0$ to be the formal sum $\sum_{i\in L}{d_i\cdot i
\over D}$.  If a 3-manifold $M^3$ is represented by a
framed link $L$ in $S^3$, then the $SO(3)$ 3-manifold invariant of $M^3$
is $\tau(M^3)={1\over D}\cdot \langle\omega_0 {*} L\rangle\cdot ({p_{-}\over
D})^{\sigma(L)}$,  where $\sigma(L)$ is the signature of the
linking matrix of $L$, and $\langle\omega_0 {*} L\rangle$ is the link invariant of $L$
where each component of $L$ is labeled by $\omega_0$ [T].  Note that in
this normalization, $\tau(S^1\times S^2)=1, \tau(S^3)={1\over D}$.
 The invariant $\tau$ is multiplicative for disjoint unions, but
 for connected sums $\tau(M_1\# M_2)=D\cdot \tau(M_1)\tau(M_2)$.

Recall that
a TQFT 3-manifold invariant is only defined for extended oriented
closed 3-manifolds, but as pointed out in [A] there is a preferred
framing for each oriented closed 3-manifold; therefore, the
formula above should be thought as for an extended 3-manifold with the preferred
framing determined by the framed link $L$.
The same
3-manifold invariant can also be defined using the representations of the mapping
class groups.  The subtlety in framing is reflected in the fact that the
TQFT representations of the mapping class groups are only
projective representations.  It is known that ${p_-\over D}$
is a root of unity of finite order [BK], so we can write
${p_-\over D}=e^{\pi i c\over 4}$ for a rational number $c$,
which is called the central charge of the TQFT (well-defined modulo 8).
Framing changes
 lead to powers of $\kappa=e^{\pi i c\over 4}$.  So up
to powers of $\kappa$, the same $\SO(3)$
invariant of 3-manifolds can be obtained as follows:  suppose an
oriented 3-manifold $M^3$ is given by gluing together two genus=$g$
oriented handlebodies
$H_g$ by a self-diffeomorphism $f$ of $\Sigma_g=\partial H_g$ (note that the
handlebody $H_g$ determines a vector $v_0$ in the TQFT vector
space $V_r(\Sigma_g)$ of $\Sigma_g$ up to a power of $\kappa$);
then the $SO(3)$ 3-manifold invariant $\tau(M^3)$ is, up to a power of $\kappa$,
the inner product of $v_0$
with $\rho_{SO(3)}(f)(v_0)$.
Theorem 3 has a direct corollary concerning the
Witten-Reshitikhin-Turaev $\SO(3)$ invariants of 3-manifolds.

\medskip
\noindent {\bf Theorem 4.}  The set of the norms of the
$\SO(3)$ invariants at $A=ie^{2\pi i
\over 4r}$  of all connected 3-manifolds is dense in $[0,\infty)$ for
primes $r\ge 5$.

\medskip

\proof
Given any complex number $z$, note that $D>1$ so we can find a $g>1$ so that $z'={z\over
  D^{g-1}}$ satisfies $\vert z'\vert < 1$.  Then arrange $z'$ as the
  $(1,1)$ entry of a unitary matrix $U_z$, and
  the vector $v_0$ determined by $H_g$ as the first basis vector of an orthonormal
  basis of $V_r(\Sigma)$.  By Theorem~3,
  $U_z$ can be approximated by a sequences of
  unitary matrices associated to diffeomorphisms $f_i$ of $\Sigma_g$ up to a phase.
  Each $f_i$ determines a 3-manifold $M_{f_i}$ by gluing two
  copies of $H_g$.
  It follows that the 3-manifold invariant of $M_{f_i}$ is the (1,1)-entry
  of $\rho_{SO(3)}(f_i)$ up to a phase, hence approximates $z'$ up to a phase.
  By connecting sums $M_{f_i}$ with $(g-1)$ copies of $S^1\times S^2$, we
  approximate $z$ using the invariants of the 3-manifolds
  $M_{f_i}\# (g-1) (S^1\times S^2)$.  \qed

\medskip

{\it Remarks:} (1): A similar density result of the $\SU(2)$
invariants for $r=1\; mod \; 4$ can be deduced from the $\SO(3)$ case and the
decomposition formula.

(2): This result does not follow from Funar's result [F].  The
infinite mapping classes are in the image of the braid groups, and
can be extended to the handlebody groups of $H_g$.  Therefore, all resulting
3-manifolds have the same invariant up to powers of $\kappa$.

Finally, we make the following:

{\bf Conjecture:} The set of the
$\SO(3)$ invariants at $A=ie^{2\pi i
\over 4r}$  of all connected 3-manifolds is dense in the complex plane for
primes $r\geq 5$.

\bigskip\goodbreak
\centerline{\bf References}
\medskip

\item{[A]}M. Atiyah: On framings of 3-manifolds, {\it Topology}
{\bf 29} (1990), 1-7.
\item{[At]}J. H. Conway, R. T. Curtis, S. P. Norton, R. A. Parker, R. A.
    Wilson: Atlas of Finite Groups, Clarendon Press, Oxford, 1985.
\item{[Bo]}N.~Bourbaki: Groupes et alg\`ebres de Lie, Chap. 4--6 (1968),
    Chap. 7--8 (1975), Hermann, Paris.
\item{[BK]}B. ~Bakalov, A. Kirillov: Lectures on tensor categories
 and modular functors, University lecture series vol. 21, AMS, 2000.
\item{[EM]}S.~Eilenberg, S.~Mac~Lane:
    Cohomology theory in abstract groups. II. Group extensions with a  non-Abelian kernel.
    {\it Ann. of Math.} {\bf 48} (1947), 326--341.
\item{[F]}L. Funar, On the TQFT
    representations of the mapping class groups.
    Pacific J. Math. 188 (1999), no. 2, 251--274.
\item{[FKLW]}M. Freedman, A. Kitaev, M. Larsen, Z. Wang:
    Topological quantum computation, Bull. AMS 40 (2003), 31-38.
\item{[FLW1]}M. Freedman, M. Larsen, Z. Wang:  A modular functor
    which is universal for quantum computation. {\it Comm. Math. Phys.}
    {\bf 227}  (2002), 605-622, arXiv: quant-ph/0001108.
\item{[FLW2]}M. Freedman, M. Larsen, Z. Wang: The two-eigenvalue problem
    and density of Jones representation of braid groups.
    {\it Comm. Math. Phys.} {\bf 228} (2002), 177-199, arXiv: math.GT/0103200.
\item{[FK]}C. Frohman, J. Kania-Bartoszynska: SO(3) topological
quantum field theory, {\it Comm. Anal. Geom.}{\bf 4} (1996), no. 4, 589--679.
\item{[FNWW]}M. Freedman, C. Nayak, K. Walker, Z.Wang: Picture TQFTs,
    in preparation.
\item{[G]}P. Gilmer: On the Witten-Reshetikhin-Turaev representations of mapping class groups.
    {\it Proc. Amer. Math. Soc.} {\bf 127} (1999), no. 8, 2483--2488.
\item{[Ge]}P. Gerardin, Weil representations associated to finite fields.
    {\it J. Algebra} {\bf 46} (1977), no. 1, 54-101.
\item{[J]}L. Jeffrey: Chern-Simons-Witten invariants of lens spaces and torus bundles,
    and the semiclassical approximation.
    {\it Comm. Math. Phys.} {\bf 147} (1992), no. 3, 563--604.
\item{[KL]}L. Kauffman, S. Lins:
    Temperley-Lieb recoupling theory and invariants of 3-manifolds.
    Ann. Math. Studies, vol 134, Princeton University Press, 1994.
\item{[MP]}W. Mckay, J. Patera:
    Tables of dimensions, indices, and branching rules for representations of simple Lie algebras,
    Lecture notes in pure and applied math., vol 69.
\item{[R]}J. Roberts:
    Irreducibility of some quantum representations of mapping class groups.
    {\it J. Knot Theory Ramifications} {\bf 10} (2001), no. 5, 763--767.
\item{[RT]}N. Reshetikhin, V. Turaev: Invariants of 3-manifolds
via link polynomials and quantum groups. {\it Invent. Math.}{\bf
103} (1991), 547-598.
\item{[Sp]}T.~A.~Springer: Characters of Special Groups, in
    {\it Seminar on Algebraic Groups and Related Finite Groups},
    Lecture Notes in Mathematics 131,
    Springer-Verlag, Berlin, 1970.
\item{[T]}V. Turaev,
    {\it Quantum invariants of knots and 3-manifolds},
    de Gruyter Studies in Math., vol 18, 1994.
\item{[Wi]}E. Witten:Quantum field theory and the Jones polynomial.
    {\it Comm. Math. Phys.} {\bf 121} (1989), no. 3, 351--399.
\item{[Wr1]}G. Wright:
    The Reshetikhin-Turaev representation of the mapping class group.
    {\it J. Knot Theory Ramifications} {\bf 3} (1994), no. 4, 547--574.
\item{[Wr2]}G. Wright:
    The Reshetikhin-Turaev representation of the mapping class group at the sixth root of unity.
    {\it J. Knot Theory Ramifications} {\bf 5} (1996), no. 5, 721--739.

\end